\documentclass[11pt]{article}
\usepackage[T1]{fontenc}
\usepackage{lmodern}
\usepackage{amsmath,amssymb,amsthm,mathtools}
\usepackage[margin=1in]{geometry}
\usepackage{booktabs,array,microtype}
\usepackage[colorlinks=true,linkcolor=blue,citecolor=blue,urlcolor=blue]{hyperref}

\newtheorem{theorem}{Theorem}[section]
\newtheorem{lemma}[theorem]{Lemma}
\newtheorem{proposition}[theorem]{Proposition}
\newtheorem{corollary}[theorem]{Corollary}

\newcommand{\Z}{\mathbb Z}
\newcommand{\F}{\mathbb F}
\newcommand{\vv}{v_3}
\newcommand{\ind}{\mathbf1}

\title{Exact $3$-Adic Periods of Symmetric Dellac Numbers\\
and Their Polynomial Extensions}
\author{
Tony Newton\\[-1mm]
{\small Newton Astro Labs, London, UK}\\
{\small\texttt{tony.newton79@gmail.com}}
\and
Jiang Zeng\\[-1mm]
{\small College of Mathematics and Physics, Wenzhou University, Wenzhou 325035, China}\\
{\small Universit\'e Claude Bernard Lyon 1, CNRS,  ICJ UMR5208,
 69622 Villeurbanne, France}\\
{\small\texttt{zeng@math.univ-lyon1.fr}}
}
\date{}

\begin{document}
\maketitle

\begin{abstract}
\emergencystretch=2em
The odd and even symmetric Dellac numbers are the specializations
$a_n=f_n(0)$ and $e_n=f_n(1)$ of the Randrianarivony--Zeng family
$f_n(s)=2^{-n}D_n(s)$.  Using a finite Stieltjes-fraction truncation, we establish
the sharp congruence $v_3(a_n)=1$ for $n\ge2$ and obtain
explicit congruences modulo $81$ for both sequences. We also establish
universal eventual periodicity modulo arbitrary positive integers.
For every integer $s\ge0$ and $k\ge1$, we give and analyse formulas for the \emph{least}
eventual period of $f_n(s)$ modulo $3^k$ and its \emph{least starting
index}.  The period is governed by five linear factors in $s$;
the onset follows a nine-class formula modulo $9$.  In particular,
for $k\ge3$ the even symmetric Dellac period is three times the
odd period and begins one index earlier.  The proofs combine a finite
Stieltjes-fraction certificate, a graded $3$-adic pole expansion,
and an explicit finite top-degree certificate over $\F_3$.
Finite integer-arithmetic checks supplement the proofs.  We also record the
connection with the type-$B$ homogenized Linial arrangement.
\end{abstract}

\section{Introduction and main results}

Congruences of combinatorial counting sequences can display regularity long
before their ordinary generating functions have a simple closed form.  For the
symmetric Dellac numbers, the practical question is not merely whether the
terms eventually repeat modulo a fixed power of $3$, but exactly how long the
repetition takes and when it first begins.  The underlying polynomial family
permits both a continued-fraction approach to divisibility and a separate
factorial-series approach to precise periods. The continued-fraction approach yields a short proof of
the observed valuation law, while the factorial-series
approach leads to exact periods and starting indices
for the entire polynomial family. 
Our principal results give explicit formulas for the exact
eventual periods and their least starting indices;
 they are distinct from the previously known
combinatorial models and general modular-congruence methods.
The polynomials $D_n(s)$ of Randrianarivony and Zeng~\cite{RZ}
have the factorial ordinary generating function
\begin{equation}\label{eq:factorial}
 F_s(t):=\sum_{n\ge0}f_n(s)t^n
 =\sum_{m\ge0}
 \frac{m!(s+1)_m\,t^m}
 {\displaystyle\prod_{j=1}^m(1+T_jt)},
 \qquad T_j=\frac{j(j+1)}2,
 \qquad f_n(s)=\frac{D_n(s)}{2^n}.
\end{equation}
Here $(s+1)_m=(s+1)\cdots(s+m)$, with the empty product $1$.
For $s\in\Z_{\ge0}$, all the coefficients $f_n(s)$ are integers.
Bigeni and Feigin~\cite{BF} identified $f_n(0)$ and $f_n(1)$
with the counts of odd and even symmetric Dellac configurations
(OEIS A098278 and A098279, respectively).  Write
$a_n=f_n(0)$ and $e_n=f_n(1)$.  Their first values are
\[
 1,1,3,21,267,\ldots,
 \qquad 1,2,10,98,1594,\ldots.
\]
In an observation recorded in OEIS in 2021, F.~Chapoton noted
that $a_n$ appeared to have valuation exactly one at the prime~$3$
for every $n\ge2$~\cite{OEIS}.  This congruence follows from a short truncation of the
Stieltjes continued fraction. We present the corresponding
finite certificate in Section~\ref{sec:finite}, before
extending the analysis to the full polynomial family.

The underlying objects and methods have an extensive history.
Randrianarivony and Zeng~\cite{RZ} introduced the polynomial family
and its recurrences; Bigeni and Feigin~\cite{BF} supplied the
symmetric Dellac interpretations of its two adjacent specializations.
Barsky and Dumont~\cite[Th\'eor\`eme~1]{BD} obtained factorial-product
expansions and prime-power congruences for Genocchi numbers of the
second kind, while Lazar and Wachs~\cite[Sec.~2.3]{LW} described
related factorial expressions in their study of the homogenized
Linial arrangement. The factors in those classical formulas involve
$j(j+1)$, whereas the normalization in~\eqref{eq:factorial} uses
$j(j+1)/2$.

Continued-fraction congruence arguments appear in Flajolet~\cite{Flajolet};
Deb and Sokal~\cite{DebSokal} obtain Stieltjes- and Thron-type
continued fractions for multivariate polynomials enumerating
$D$-permutations and related Genocchi families. Their emphasis is
combinatorial enumeration and continued-fraction refinement, rather
than determination of the least eventual periods and their starting
indices for the present family. Krattenthaler and
M\"uller~\cite{KM3,KMp} developed methods for prime-power congruences
of sequences satisfying algebraic differential equations; their
framework provides relevant methodological context but is not invoked
in the proofs below. Han and Zeng's $q$-Genocchi divisibility
results~\cite{HZ} concern another specialization. Against this
background, the claims of the present paper are the explicit minimal
$3$-power periods and least onsets for all $s\ge0$, together with
the sharpened odd and even specializations. This comparison identifies
the scope of the results without claiming that it exhausts the
literature.

Write $T_k(s)$ for the least eventual period of $f_n(s)$ modulo
$3^k$, and $N_k(s)$ for the least nonnegative integer such that
\begin{equation}\label{eq:least-onset}
 f_{n+T_k(s)}(s)\equiv f_n(s)\pmod{3^k}
 \quad\text{for every }n\ge N_k(s).
\end{equation}
For $s\ge0$, define
\begin{equation}\label{eq:tau-rho}
 \begin{split}
 \mathcal P(s)&=(s+1)(s+2)(s+3)(s+4)(s+6),\\
 \tau(s)&=1+\vv(\mathcal P(s)),\\
 \rho(s)&=\begin{cases}
  \vv(s+3),&s\equiv0\pmod3,\\
  0,&s\equiv1\pmod3,\\
  \vv(s+1),&s\equiv2\pmod3.
 \end{cases}
 \end{split}
\end{equation}
All these valuations are finite, since $s\ge0$.

\begin{theorem}[Exact eventual periods]\label{thm:period}
For every $s\in\Z_{\ge0}$ and $k\ge1$,
\begin{equation}\label{eq:period}
\displaystyle
 T_k(s)=\begin{cases}
  1,&1\le k\le\rho(s),\\
  2\cdot3^{\max\{0,k-\tau(s)\}},&k>\rho(s).
 \end{cases}
\end{equation}
\end{theorem}

For compactness, put
\begin{equation}\label{eq:odd-onset}
 \mathcal N(j)=\begin{cases}
 2,&j=1,\\
 j,&j\ge2,\ j\equiv1\pmod3,\\
 j+1,&j\ge2,\ j\not\equiv1\pmod3.
 \end{cases}
\end{equation}

\begin{theorem}[Exact least starting indices]\label{thm:onset}
For every $s\in\Z_{\ge0}$ and $k\ge1$, the least index
in~\eqref{eq:least-onset} is given by the following table.
In the last two rows, $r$ is as specified there.
\begin{equation}\label{eq:onset-table}
\begin{array}{c|l}
 s\pmod9 & N_k(s)\\ \hline
0 & \mathcal N(k)\\
1 & 0\ (k=1);\quad k-\ind_{k\equiv1\ (3)}\ (k\ge2)\\
2 & 1\ (k=1);\quad k-\ind_{k\equiv1\ (3)}\ (k\ge2)\\
3 & k+1-\ind_{k\equiv0\ (3)}\\
4 & 0\ (k=1);\quad k\ (k\ge2)\\
5 & k\\
6 & k+1\\
7 & 0\ (k\le r);\quad\mathcal N(k-r)\ (k>r),\quad r=\vv(s+2)\\
8 & 1\ (k\le r);\quad\mathcal N(k-r)\ (k>r),\quad r=\vv(s+1)
\end{array}
\end{equation}
\end{theorem}

\begin{corollary}[Symmetric Dellac specializations]
For $k\ge3$, the odd specialization $s=0$ has least eventual
period $2\cdot3^{k-3}$, and for $k\ge2$ the even specialization
$s=1$ has least eventual period $2\cdot3^{k-2}$.
Their respective least starting indices are
$\mathcal N(k)$ and $k-\ind_{k\equiv1\ (3)}$ for $k\ge2$.
In particular, for every $k\ge3$,
\[
 T_k(1)=3T_k(0),\qquad N_k(1)=N_k(0)-1.
\]
\end{corollary}

\section{Finite Stieltjes fractions and low-power congruences}
\label{sec:finite}

The polynomial family has the Stieltjes fraction
\begin{equation}\label{eq:general-S}
 F_s(t)=\cfrac{1}{1-\cfrac{\alpha_1(s)t}{
                 1-\cfrac{\alpha_2(s)t}{
                 1-\cfrac{\alpha_3(s)t}{\ddots}}}},
 \quad
 \alpha_{2r-1}(s)=r(2r-1+s),\quad
 \alpha_{2r}(s)=r(2r+s)
 \quad (r\ge1),
\end{equation}
as follows from the Randrianarivony--Zeng framework~\cite{RZ}.
Thus its odd specialization has coefficients
$1,2,6,8,15,18,\ldots$, while its even specialization has
$2,3,8,10,18,21,\ldots$.
Write $R_h^{(s)}=P_h^{(s)}/Q_h^{(s)}$ for the finite fraction
through $\alpha_h(s)$, with $R_0^{(s)}=1$.
The continuants satisfy
\begin{equation}\label{eq:continuants}
 P_0^{(s)}=Q_0^{(s)}=P_1^{(s)}=1,\quad
 Q_1^{(s)}=1-\alpha_1(s)t,\quad
 X_h^{(s)}=X_{h-1}^{(s)}-\alpha_h(s)tX_{h-2}^{(s)}
 \quad(h\ge2)
\end{equation}
for $X=P,Q$.

\begin{lemma}[Product truncation]\label{lem:product}
Let $C_h(s)=\prod_{j=1}^h\alpha_j(s)$.  Then
\begin{equation}\label{eq:product-difference}
 R_h^{(s)}-R_{h-1}^{(s)}
 =\frac{C_h(s)t^h}{Q_h^{(s)}Q_{h-1}^{(s)}},
 \qquad
 F_s(t)-R_{h-1}^{(s)}(t)\in C_h(s)t^h\Z[[t]].
\end{equation}
In particular, if $M\mid C_h(s)$, then
$F_s(t)\equiv R_{h-1}^{(s)}(t)\pmod M$.
\end{lemma}

\begin{proof}
The continuant recurrence gives
$P_h^{(s)}Q_{h-1}^{(s)}-P_{h-1}^{(s)}Q_h^{(s)}
=C_h(s)t^h$ by induction.  Each $Q_h^{(s)}$ has
constant coefficient~$1$, so is a unit in $\Z[[t]]$.
The convergents converge $t$-adically to $F_s$, and
the differences for $j\ge h$ are all divisible by
$C_h(s)t^h$.
\end{proof}

\begin{proposition}[Universal eventual modular periodicity]
\label{prop:allmoduli}
For every integer $s\ge0$ and positive integer $M$,
the sequence $\{f_n(s)\bmod M\}_{n\ge0}$ is
eventually periodic.
\end{proposition}

\begin{proof}
Since $\alpha_{2M}(s)=M(2M+s)$, the Stieltjes
fraction terminates modulo $M$ at that coefficient.
Equivalently, Lemma~\ref{lem:product} gives a finite
rational series $F_s(t)\equiv P(t)/Q(t)\pmod M$,
where $Q(0)=1$.  Its coefficients eventually satisfy
a finite linear recurrence over $\Z/M\Z$.
The recurrence has finitely many states, so some state
eventually repeats.
\end{proof}

The following proof provides a short finite-fraction
certificate for the sharp congruence modulo $9$.
The key observation is that $\alpha_6(0)=18$ vanishes
modulo $9$, allowing the Stieltjes fraction to be
truncated at this coefficient.

\begin{theorem}[Sharp congruence for odd symmetric Dellac numbers]
\label{thm:odd9}
For the odd specialization,
\begin{equation}\label{eq:odd9}
 a_0\equiv a_1\equiv1,\qquad a_2\equiv3,\qquad
 a_n\equiv3(-1)^{n-3}\quad(n\ge3)\pmod9.
\end{equation}
Consequently $v_3(a_n)=1$ for every $n\ge2$.
\end{theorem}

\begin{proof}
By Lemma~\ref{lem:product} (or directly by
$\alpha_6(0)\equiv0\pmod9$),
\[
 F_0(t)\equiv R_5^{(0)}(t)
 =\frac{N(t)}{D(t)}
 =\frac{1-31t+136t^2}{1-32t+165t^2-90t^3}
 \pmod9.
\]
Put $H(t)=1+2t+4t^2+6t^3$.
The required polynomial identity is
\begin{equation}\label{eq:odd9-certificate}
 N(t)(1+t)-D(t)H(t)
 =18t^3(30t^3-35t^2-16t+1).
\end{equation}
As $D(0)=1$, the certificate gives
\[
 F_0(t)\equiv\frac{H(t)}{1+t}
 =1+t+3t^2+\frac{3t^3}{1+t}\pmod9,
\]
whose coefficients have the asserted residues.
\end{proof}

\begin{theorem}[Explicit congruences modulo $81$]
\label{thm:both81}
The odd and even specializations satisfy
\begin{align}
 a_n&\equiv(-1)^n(51-27n)\pmod{81}
 &&(n\ge4),\label{eq:odd81}\\
 e_n&\equiv(-1)^n(10-9n)\pmod{81}
 &&(n\ge3).\label{eq:even81}
\end{align}
Moreover, $e_n\equiv(-1)^n\pmod3$ for every $n\ge0$
and $e_n\equiv(-1)^n\pmod9$ for $n\ge2$.
In particular, every $e_n$ is prime to~$3$.
\end{theorem}

\begin{proof}
For $s=0$, the product of the first six Stieltjes
coefficients is $25920$, divisible by~$81$.
Thus $F_0\equiv R_5^{(0)}\pmod{81}$, with $N,D$ as
in the preceding proof.  For
$H_0=1+3t+6t^2+28t^3-12t^4-9t^5$,
direct multiplication yields
\begin{equation}\label{eq:odd81-certificate}
 N(1+t)^2-DH_0
 =81t^4(4-55t+52t^2+5t^3-10t^4).
\end{equation}
Hence $F_0\equiv H_0/(1+t)^2\pmod{81}$.
The coefficients of $(1+t)^{-2}$ are
$(-1)^n(n+1)$; extraction gives~\eqref{eq:odd81}.

For $s=1$, the first six Stieltjes coefficients have
product $181440$, again divisible by~$81$, and
\[
 R_5^{(1)}(t)
 =\frac{P(t)}{Q(t)}
 =\frac{1-39t+228t^2}{1-41t+300t^2-288t^3}.
\]
With $H_1=1+4t+15t^2+39t^3+18t^4$ we have
\begin{equation}\label{eq:even81-certificate}
 P(1+t)^2-QH_1
 =81t^3(1-19t-82t^2+72t^3+64t^4).
\end{equation}
Thus $F_1\equiv H_1/(1+t)^2\pmod{81}$, giving
\eqref{eq:even81}.  Finally, in the factorial series
\eqref{eq:factorial}, the terms $m\ge2$ vanish modulo
$3$ for $s=1$; hence
$F_1(t)\equiv1+2t/(1+t)=(1+t)^{-1}\pmod3$.
Modulo~$9$, terms $m\ge3$ vanish and
\[
 F_1(t)\equiv 1+\frac{2t}{1+t}
 +\frac{12t^2}{(1+t)(1+3t)}
 \equiv\frac{1+3t+3t^2}{1+t}\pmod9.
\]
This proves the remaining assertions.
\end{proof}

\section{The graded pole expansion}

Put
\begin{equation}\label{eq:Pq}
 P_m(s)=m!(s+1)_m=(m!)^2\binom{s+m}{m},
 \qquad q_m=\left\lfloor\frac{m+2}{3}\right\rfloor.
\end{equation}
Exactly $q_m$ factors in the denominator of the $m$th
summand of~\eqref{eq:factorial} are congruent to $1+t$
modulo $3$.  In fact,
\begin{equation}\label{eq:triangular}
 T_{3a+1}=1+9\frac{a(a+1)}2,
 \qquad T_j\equiv0\pmod3\quad(j\equiv0,2\pmod3).
\end{equation}
When $3\mid T_j$, the factor $(1+T_jt)^{-1}$ is a polynomial
modulo $3^k$.  For $j=3a+1$, expand
\begin{equation}\label{eq:unit-pole}
 \frac1{1+T_jt}
 =\frac1{1+t}\sum_{h\ge0}
 \left(-\frac{9a(a+1)t}{2(1+t)}\right)^h.
\end{equation}
Thus each additional pole at $1+t=0$ costs at least two
powers of $3$.  All expansions terminate modulo a fixed $3^k$.

\begin{lemma}[Eventual signed binomial expansion]\label{lem:binomial}
For fixed $s\ge0$ and $k\ge1$, there is an integer $J$,
an index $n_0$, and integers $c_{k,0}(s),\ldots,c_{k,J}(s)$
such that, for all $n\ge n_0$,
\begin{equation}\label{eq:binomial}
 f_n(s)\equiv(-1)^n\sum_{j=0}^J
 c_{k,j}(s)\binom nj\pmod{3^k}.
\end{equation}
Their coefficients can be chosen so that
\begin{equation}\label{eq:graded-bound}
 \vv(c_{k,j}(s))\ge \tau(s)+2j-2
 \quad(j\ge1),
\end{equation}
with the usual convention that the inequality means divisibility
when a coefficient vanishes modulo $3^k$.
Moreover, the constant coefficient is divisible by
$3^{\rho(s)}$.
\end{lemma}

\begin{proof}
For all $s\ge0$, $\vv(P_m(s))\ge2\vv(m!)\ge2q_m-2$.
It follows from~\eqref{eq:unit-pole} that, without further
refinement, a pole of order $q$ is accompanied by at least
$3^{2q-2}$.  This gives~\eqref{eq:graded-bound} when
$s\equiv1,4\pmod9$, for which $\tau(s)=2$.

When $s\equiv2\pmod3$, the product $(s+1)_m$ contains
one more multiple of $3$ than the crude bound for
$m\equiv1,2\pmod3$; hence $\vv(P_m(s))\ge2q_m-1$.
When $s\equiv0\pmod3$, pair $m=3a+1$ with $m+1=3a+2$:
\begin{equation}\label{eq:pair-general}
 F_{s,m}(t)+F_{s,m+1}(t)
 =\frac{P_m(s)t^m\bigl((1+t)+3d_a(s)t\bigr)}
 {\prod_{j=1}^{m+1}(1+T_jt)},
\end{equation}
where
\[
 d_a(s)=\frac{T_{m+1}+(m+1)(s+m+1)-1}{3}\in\Z.
\]
The $(1+t)$ part cancels one pole, while the second part
gains a factor $3$.  Terms with $m\equiv0\pmod3$ already
have sufficient factorial valuation.  Therefore for
$s\equiv0,2\pmod3$, every pole of order $q\ge1$
comes with at least $3^{2q-1}$.  This establishes the
claimed estimate for $s\equiv0,2\pmod9$.

It remains to refine the other five classes.  For
$s\equiv7\pmod9$, put $u=s+2$, $r=\vv(u)\ge2$.
For $m\ge2$, $P_m(s)=uQ_m$ with
\[
 Q_m=m!(u-1)\prod_{j=1}^{m-2}(u+j),\qquad
 \vv(Q_m)\ge\vv(m!)+\vv((m-2)!)\ge2q_m-3.
\]
Thus poles of order $q\ge2$ carry at least
$3^{r+2q-3}$, as required because $\tau(s)=r+1$.

For $s\equiv8\pmod9$, put $u=s+1$, $r=\vv(u)\ge2$.
Now $P_m(s)=uQ_m$ for $m\ge1$, with
\[
 Q_m=m!\prod_{j=1}^{m-1}(u+j),\qquad
 \vv(Q_m)\ge\vv(m!)+\vv((m-1)!)\ge2q_m-2.
\]
Consequently poles of order $q\ge2$ carry
$3^{r+2q-2}$, and $\tau(s)=r+2$.

For $s\equiv5\pmod9$, take $u=s+4$, $r=\vv(u)\ge2$.
All $m\ge4$ have
\[
 P_m(s)=uQ_m,\qquad
 Q_m=m!(u-3)(u-2)(u-1)
       \prod_{j=1}^{m-4}(u+j).
\]
Here $\vv(u-3)=1$, so
$\vv(Q_m)\ge\vv(m!)+1+\vv((m-4)!)\ge2q_m-2$.
The terms $m\le3$ have at most one unit pole.
This gives the bound $3^{r+2q-2}$ for $q\ge2$.

For $s\equiv6\pmod9$, set $u=s+3$, $r=\vv(u)\ge2$.
First, the exact identity
\begin{equation}\label{eq:minus3}
 F_{s,1}(t)+F_{s,2}(t)
 =\frac{(u-2)t}{1+3t}
  +\frac{2u(u-2)t^2}{(1+t)(1+3t)}
\end{equation}
shows that the low terms have no higher-order unit poles
without a factor $u$.  For $m\ge3$,
\[
 P_m(s)=uQ_m,\quad
 Q_m=m!(u-2)(u-1)\prod_{j=1}^{m-3}(u+j),
 \quad\vv(Q_m)\ge\vv(m!)+\vv((m-3)!).
\]
For $m\equiv0\pmod3$ this suffices directly.
For $m=3a+1,3a+2$, the bare lower bound is only
$r+2q_m-3$.  In the pair~\eqref{eq:pair-general}, the
$(1+t)$ part removes one of the $q_m$ unit poles, while the
$3d_a(s)t$ part gains one power of $3$ without changing
that pole count.  Both contributions therefore satisfy the
required bound $r+2q-2$ for their \emph{actual} pole order
$q$.  Consequently all poles with $q\ge2$ cost at least
$3^{r+2q-2}$.

Finally, for $s\equiv3\pmod9$, put $u=s+6$,
$r=\vv(u)\ge2$.  The exceptional terms $m=4,5$ satisfy
the exact cancellation
\begin{align}\label{eq:minus6}
 F_{s,4}(t)+F_{s,5}(t)
 ={}&\frac{P_4(s)t^4}
 {(1+t)(1+3t)(1+6t)(1+15t)}\notag\\
 &+\frac{5uP_4(s)t^5}
 {\prod_{j=1}^{5}(1+T_jt)}.
\end{align}
Since $\vv(P_4(s))=2$, its only double-pole contribution
is divisible by $3^{r+2}$.  For $m\ge6$,
\[
 P_m(s)=uQ_m,\quad
 Q_m=m!(u-5)(u-4)(u-3)(u-2)(u-1)
       \prod_{j=1}^{m-6}(u+j),
\]
with $\vv(Q_m)\ge\vv(m!)+1+\vv((m-6)!)$.
For $m=3a+1,3a+2$ beginning at $m=7$, the bare
factorial bound is again one power of $3$ short.  Pairing
these consecutive terms by~\eqref{eq:pair-general}
removes an exact $(1+t)$ factor or supplies an extra factor
$3$, as in the $s\equiv6\pmod9$ case.  Thus all double
and higher poles cost at least $3^{r+2q-2}$.


\medskip
\noindent\textit{Uniform verification of the exceptional
valuation bounds.}
We make explicit the factorial estimates used above.

Let $u=s+d$, $r=\vv(u)\ge2$, where
$d\in\{1,2,3,4,6\}$ corresponds to the five
exceptional residue classes.
Since
\[
 \binom{u+h}{h}\in\Z
 \qquad(u\ge0,\ h\ge0),
\]
we have
\begin{equation}\label{eq:factorial-uniform}
 \vv\left(\prod_{j=1}^{h}(u+j)\right)
 \ge \vv(h!)\ge\left\lfloor\frac h3\right\rfloor.
\end{equation}

For $m\ge d$, the factorization
\[
 P_m(s)=m!\,u
 \left(\prod_{j=1}^{d-1}(u-j)\right)
 \left(\prod_{j=1}^{m-d}(u+j)\right)
\]
gives
\begin{equation}\label{eq:uniform-P-bound}
 \vv(P_m(s))
 \ge r+\vv(m!)+\vv((m-d)!)+b_d,
\end{equation}
where
\[
 b_d=\sum_{j=1}^{d-1}\vv(j).
\]
Indeed, since $9\mid u$ and $1\le j<d\le6$,
we have $\vv(u-j)=\vv(j)$.

For $d\in\{1,2,3,4,6\}$, the required
valuation of a pole of order $q\ge2$ is
\[
 \tau(s)+2q-4=r+2q-c_d,
\]
where
\[
 c_2=3,\qquad c_1=c_3=c_4=c_6=2.
\]

Using
\[
 \vv(m!)\ge\left\lfloor\frac m3\right\rfloor,
 \qquad
 \vv((m-d)!)\ge
 \left\lfloor\frac{m-d}{3}\right\rfloor,
\]
define
\[
 \varepsilon_d(m)=
 \left\lfloor\frac m3\right\rfloor+
 \left\lfloor\frac{m-d}{3}\right\rfloor+
 b_d-(2q_m-c_d).
\]
A direct calculation gives
\[
\begin{array}{c|c|c|ccc}
d&s\pmod9&b_d&
 m\equiv0&m\equiv1&m\equiv2\\ \hline
2&7&0&2&0&1\\
1&8&0&1&0&0\\
4&5&1&1&0&0\\
3&6&0&1&-1&-1\\
6&3&1&1&-1&-1
\end{array}
\]
where the last three columns give
$\varepsilon_d(m)$ according to $m\bmod3$.

Thus the required valuation holds termwise
for $d=1,2,4$.  For $d=3,6$, it also
holds termwise when $m\equiv0\pmod3$.

It remains to consider $d=3,6$ and
$m\equiv1,2\pmod3$.  Pair $m=3a+1$
with $m+1=3a+2$, so that
$q_m=q_{m+1}$.  The exact identity
\eqref{eq:pair-general} reads
\[
 F_{s,m}(t)+F_{s,m+1}(t)
 =
 \frac{P_m(s)t^m
       \bigl((1+t)+3d_a(s)t\bigr)}
      {\prod_{j=1}^{m+1}(1+T_jt)}.
\]
The term containing $(1+t)$ cancels
one unit pole, lowering its order by one.
The other term gains an additional
factor $3$.  In either case, the
one-power deficit in the table is removed.

Finally, each additional unit-pole
correction in~\eqref{eq:unit-pole}
raises the pole order by one while
contributing at least two additional
powers of $3$.  Hence the required
valuation bound remains valid for
every higher pole order.

The exceptional low-index terms are
already covered by the exact identities
\eqref{eq:minus3} and~\eqref{eq:minus6}.
Consequently, every pole of order
$q\ge2$ has coefficient divisible by
\[
 3^{\tau(s)+2q-4}.
\]

To finish, a summand $t^h/(1+t)^q$ has for large $n$
the coefficient
$(-1)^{n-h}\binom{n-h+q-1}{q-1}$.
The latter polynomial has an integral expansion in the basis
$\binom nj$, of degrees at most $q-1$.
Thus the coefficient of $\binom nj$ receives contributions
only from poles of order $q\ge j+1$.
The preceding estimates imply~\eqref{eq:graded-bound}.
The statement about the constant coefficient follows from
the low-power computations in the next lemma.
\end{proof}

\section{Sharp valuations of the constant and linear coefficients}

\begin{lemma}[The constant coefficient]\label{lem:constant}
For $k>\rho(s)$, the constant $c_{k,0}(s)$ in
Lemma~\ref{lem:binomial} has valuation exactly $\rho(s)$.
For $k\le\rho(s)$, the entire sequence is eventually zero
modulo $3^k$.
\end{lemma}

\begin{proof}
Modulo $3$, the factorial expansion gives
\begin{equation}\label{eq:mod3}
 F_s(t)\equiv\begin{cases}
 1+t,&s\equiv0\pmod3,\\
 (1+t)^{-1},&s\equiv1\pmod3,\\
 1,&s\equiv2\pmod3.
 \end{cases}
\end{equation}
For $n$ sufficiently large, expansion through $m=2$
modulo $9$ yields
\begin{equation}\label{eq:constant9}
 f_n(s)\equiv(-1)^n(s+1)(8s+15)\pmod9.
\end{equation}
This is a unit for $s\equiv1,4,7\pmod9$, and has
valuation exactly $1$ for $s\equiv0,2,3,5\pmod9$.

For $s\equiv8\pmod9$, write $u=s+1$, $r=\vv(u)$.
Since $F_s(t)=1+uG_u(t)$, with $G_u(t)$ obtained by
dividing all positive-index summands by $u$, one has
\[
 G_u(t)\equiv\frac{t+2t^2}{1+t}\pmod3.
\]
Its eventual signed constant is $-1+2=1\pmod3$.
Hence the constant coefficient of $F_s$ has valuation $r$.

For $s\equiv6\pmod9$, write $u=s+3$, $r=\vv(u)$.
By~\eqref{eq:minus3}, the part without $u$ has
denominator $1+3t$ and eventually vanishes modulo $3^k$.
After dividing the other terms by $u$, their eventual
signed constant is $2\pmod3$, from the second term
of~\eqref{eq:minus3}; the $m\ge3$ terms vanish modulo $3$.
Thus the constant coefficient has valuation exactly $r$.
The positive-degree coefficients have strictly larger
valuation by Lemma~\ref{lem:binomial}, completing the proof.
\end{proof}

\begin{lemma}[The linear coefficient]\label{lem:slope}
Whenever $k>\tau(s)$, the coefficient $c_{k,1}(s)$
in~\eqref{eq:binomial} has $3$-adic valuation exactly
$\tau(s)$.
\end{lemma}

\begin{proof}
The terms with $m\ge6$ vanish modulo $81$, because
$\vv(P_m(s))\ge2\vv(m!)\ge4$.
Terms with $m\le3$ contribute only simple unit poles.
Expanding the terms $m=4,5$ modulo $81$ and retaining the
double-pole part gives
\[
 \frac{P_4(s)t^4+P_5(s)t^5(1-24t)}{(1+t)^2}.
\]
Therefore the coefficient of $n$ in the eventual signed
polynomial satisfies
\begin{equation}\label{eq:slope81}
 c_{k,1}(s)\equiv P_4(s)-25P_5(s)\pmod{81}
 \qquad(k\ge4).
\end{equation}
The same calculation, reduced modulo $27$, applies
when $k=3$. Thus the congruence determines the
linear coefficient at every precision needed in
the four residue classes considered below.

Since
\[
 P_5(s)=5(s+5)P_4(s),
\]
equation~\eqref{eq:slope81} becomes
\[
 c_{k,1}(s)\equiv
 -P_4(s)(125s+624)\pmod{81}.
\]
For $s\equiv0\pmod9$, the two factors have
valuations $2$ and $1$, respectively.
For $s\equiv2\pmod9$, their valuations are
$3$ and $0$.  Thus the slope has valuation
$3$ in both cases.

For $s\equiv1,4\pmod9$, $P_4(s)$ has
valuation $2$ and $125s+624$ is a
$3$-adic unit.  The slope therefore has
valuation $2$.
These four valuations agree with $\tau(s)$.

It remains to consider
$s\equiv3,5,6,7,8\pmod9$.
Put $u=s+d$, where $d=6,4,3,2,1$,
respectively, and let
$r=\vv(u)\ge2$.

We work modulo $3^{\tau(s)+1}$.
An important simplification is that only
a few factorial summands can contribute.

Indeed,
\[
 P_{m+1}(s)=(m+1)(s+m+1)P_m(s),
\]
so the valuations $\vv(P_m(s))$ are
nondecreasing in $m$.

Directly from the factorizations of $P_m(s)$,
using $9\mid u$, we obtain
\[
\begin{array}{c|c|c|c}
s\bmod9&d&m_0&\vv(P_{m_0}(s))\\ \hline
3&6&6&r+3\\
5&4&6&r+3\\
6&3&6&r+3\\
7&2&5&r+2\\
8&1&6&r+3
\end{array}
\]
Consequently, all summands with $m\ge m_0$
vanish modulo $3^{\tau(s)+1}$.

The terms with $m\le3$ have at most one
unit pole, namely the factor $1+t$.
They therefore contribute no linear term
to the eventual signed binomial expansion.

For $m=4,5$, the original unit-pole
order is two.  Replacing the factor
$1+T_4t=1+10t$ by $1+t$ introduces
corrections divisible by $9$.
The factors with $3\mid T_j$ have
polynomial inverses congruent to $1$
modulo $3$.

By Lemma~\ref{lem:binomial}, poles of
order at least three have coefficients
divisible by $3^{\tau(s)+2}$.
Hence, at the precision considered here,
only the leading double-pole contributions
need to be retained.

Recall that the coefficient of $\binom n1$
in the eventual signed expansion of
$t^h/(1+t)^2$ is $(-1)^h$.

For $s\equiv3\pmod9$, use the exact
cancellation~\eqref{eq:minus6}.
Its first summand has only one unit pole.
The second summand gives the leading
linear contribution
\[
 -5uP_4(s).
\]

For $s\equiv6\pmod9$, apply
\eqref{eq:pair-general} with $m=4$.
The part containing $1+t$ has no
uncorrected double pole; its additional
unit-pole corrections are divisible
by $9$ and vanish at the required precision.
The other part gives
\[
 -3d_1(s)uQ_4,
 \qquad Q_4=P_4(s)/u.
\]

In the remaining classes the leading
linear contributions are
\[
 u(Q_4-Q_5),\quad uQ_4,\quad
 u(Q_4-Q_5)
\]
for $s\equiv5,7,8\pmod9$, respectively.

It remains to prove that these contributions
have exactly the asserted valuations,
uniformly in $r$.

After division by $u$, they admit the
following exact polynomial factorizations:
\begin{equation}\label{eq:uniform-slope-factorizations}
\begin{array}{c|l}
s\bmod9 & \text{leading slope divided by }u\\ \hline
3&
 -120(u-5)(u-4)(u-3)(u-2)\\[2pt]
5&
 -24(u-3)(u-2)(u-1)(5u+4)\\[2pt]
6&
 -24(u-2)(u-1)(u+1)(5u+24)\\[2pt]
7&
 24(u-1)(u+1)(u+2)\\[2pt]
8&
 -24(u+1)(u+2)(u+3)(5u+19)
\end{array}
\end{equation}

For the third row we used
\[
 3d_1(s)=5u+24.
\]

Since $9\mid u$, we have
\[
 \vv(u-3)=\vv(u+3)=1,
 \qquad \vv(5u+24)=1.
\]
Every other nonconstant factor in
\eqref{eq:uniform-slope-factorizations}
is a $3$-adic unit.

It follows that the five displayed
polynomials have valuations, respectively,
\[
 2,\quad2,\quad2,\quad1,\quad2.
\]

Multiplying by $u$, whose valuation is
$r$, shows that the leading linear
coefficient has valuation
\[
 \begin{cases}
 r+1,&s\equiv7\pmod9,\\
 r+2,&s\equiv3,5,6,8\pmod9.
 \end{cases}
\]
These values are precisely $\tau(s)$.

The higher-pole contributions and the
discarded factorial summands are divisible
by $3^{\tau(s)+1}$, so they cannot cancel
the displayed leading coefficients.
Thus
\[
 \vv(c_{k,1}(s))=\tau(s)
 \qquad(k>\tau(s)),
\]
as required.
\end{proof}

\section{Proof of the exact-period theorem}

\begin{proof}[Proof of Theorem~\ref{thm:period}]
Lemmas~\ref{lem:binomial}--\ref{lem:slope}
give~\eqref{eq:binomial} with
\begin{equation}\label{eq:key-filtration}
 \vv(c_{k,0})=\rho(s)\quad(k>\rho(s)),\qquad
 \vv(c_{k,1})=\tau(s)\quad(k>\tau(s)),\qquad
 \vv(c_{k,j})\ge\tau(s)+2j-2\quad(j\ge1).
\end{equation}
If $k\le\rho(s)$, all sufficiently late coefficients
vanish modulo $3^k$; hence the least eventual period is $1$.
If $\rho(s)<k\le\tau(s)$, the tail is a nonzero
alternating constant and has least period $2$.

Now suppose $k>\tau(s)$ and set
$d=2\cdot3^{k-\tau(s)}$.
For every $1\le i\le j$,
\[
 \vv\binom di\ge\vv(d)-\vv(i)
       =k-\tau(s)-\vv(i).
\]
By Vandermonde's identity,
\[
 \binom{n+d}{j}-\binom nj
  =\sum_{i=1}^j\binom di\binom n{j-i},
\]
and the contribution from the $j$th term of
\eqref{eq:binomial} is divisible by
$3^{k+2j-2-\vv(i)}$.  Since
$\vv(i)\le2j-2$, $d$ is an eventual period.

Any eventual period must be even, since the sequence
is a nonzero alternating constant modulo $3^{\rho(s)+1}$.
Suppose an even eventual period $e$ satisfies
$\vv(e)<k-\tau(s)$.
The $j=1$ term in the difference of the signed polynomials
has valuation exactly $\tau(s)+\vv(e)<k$.
In the Vandermonde expansion of each $j\ge2$ term, every
summand has valuation at least
$\tau(s)+\vv(e)+2j-2-\vv(i)$ for some $1\le i\le j$.
Since $2j-2-\vv(i)\ge1$ when $j\ge2$, the entire
contribution of all $j\ge2$ terms is divisible by
$3^{\tau(s)+\vv(e)+1}$.  The $j=1$ difference is
exactly $c_{k,1}e$ and has valuation
$\tau(s)+\vv(e)<k$.  Higher terms therefore cannot
cancel it modulo $3^k$.  No such $e$ is a period, and the
exhibited power-of-$3$ period is minimal.
\end{proof}

\section{A finite certificate for the starting indices}

The period proof concerns only eventual coefficients.
To locate the first index of periodicity we compute the
highest surviving polynomial degree after clearing the
unit poles.

Put
\begin{equation}\label{eq:degreeM}
 v_m(s)=\vv(P_m(s)),\qquad
 \Delta_m(s)=m-q_m-v_m(s),\qquad
 M(s)=\max_{m\ge2}\Delta_m(s).
\end{equation}
By~\eqref{eq:unit-pole}, $(1+t)^kF_s(t)$ is a polynomial
modulo $3^k$.  In a term with $m\ge2$, let $\ell$ count
the degree contributed by the factors with $3\mid T_j$,
and let $h$ count extra corrections from the factors
$T_j\equiv1\pmod9$.  If it survives modulo $3^k$, then
\[
 v_m(s)+\ell+2h\le k-1,
\]
and after multiplication by $(1+t)^k$ its polynomial degree
is at most $k+m-q_m+\ell$.
Consequently
\begin{equation}\label{eq:degree-bound}
 \deg\bigl((1+t)^kF_s(t)\bmod3^k\bigr)
 \le2k-1+M(s),
\end{equation}
apart from the harmless $m=0,1$ terms in the low-modulus
cases specified below.

\paragraph{Uniform cutoff at $m=8$.}
Since
\[
 (s+1)_m=m!\binom{s+m}{m},
\]
we have $v_m(s)\ge 2\vv(m!)$ for every
nonnegative integer $s$.  We claim that
\[
 2\vv(m!)\ge m-q_m+1 \qquad(m\ge9).
\]
This is immediate for $m=9,10,11$.  It then follows
by induction in steps of three: indeed,
\[
 (m+3)-q_{m+3}=(m-q_m)+2,
\]
whereas $\vv((m+3)!)-\vv(m!)\ge1$, since
any three consecutive integers contain a multiple of $3$.
Consequently $\Delta_m(s)\le-1$ for $m\ge9$.

For the exceptional classes we need a stronger bound.
Suppose first that $s\equiv7\pmod9$, and put
$u=s+2$, $r=\vv(u)\ge2$.  For $m\ge2$,
\[
 P_m(s)
 =m!\,u(u-1)(m-2)!
   \binom{u+m-2}{m-2}.
\]
As $u-1$ is a $3$-adic unit,
\[
 v_m(s)\ge r+\vv(m!)+\vv((m-2)!).
\]
Similarly, if $s\equiv8\pmod9$, put
$u=s+1$, $r=\vv(u)\ge2$.  Then
\[
 P_m(s)
 =m!\,u(m-1)!
   \binom{u+m-1}{m-1},
\]
and hence
\[
 v_m(s)\ge r+\vv(m!)+\vv((m-1)!).
\]

For $d=1,2$ one has
\[
 \vv(m!)+\vv((m-d)!)\ge m-q_m
 \qquad(m\ge9).
\]
Indeed, the inequality holds directly for
$m=9,10,11$, and is preserved under $m\mapsto m+3$:
each factorial valuation increases by at least one,
while $m-q_m$ increases by exactly two.
Thus in both exceptional classes
\[
 \Delta_m(s)\le-r<1-r
 \qquad(m\ge9).
\]
It follows that, in every residue class modulo $9$,
the maximum $M(s)$ is attained among $2\le m\le8$.

%
%
Define over $\F_3$ the rational series
\begin{equation}\label{eq:Gm}
 G_m(x)=\prod_{\substack{1\le j\le m\\3\mid T_j}}
       \left(1+\frac{T_j}{3}x\right)^{-1}.
\end{equation}
In particular,
\begin{equation}\label{eq:Gvalues}
 \begin{split}
 G_2&=(1+x)^{-1},\\
 G_3=G_4&=(1-x^2)^{-1},\\
 G_5&=((1+x)(1-x)^2)^{-1},\\
 G_6=G_7=G_8&=(1-x^2)^{-2}.
 \end{split}
\end{equation}
Write $\overline U_m(s)$ for the reduction modulo $3$ of
$P_m(s)/3^{v_m(s)}$, and set
\begin{equation}\label{eq:Hs-def}
 H_s(x)=\sum_{\substack{m\ge2\\\Delta_m(s)=M(s)}}
       \overline U_m(s)x^{v_m(s)}G_m(x)
       \quad\in\F_3[[x]].
\end{equation}
The top-degree coefficient of~\eqref{eq:degree-bound} can
arise only if $h=0$, $\ell=k-1-v_m$, and
$\Delta_m=M(s)$.  Hence
\begin{equation}\label{eq:top-certificate}
 3^{1-k}[t^{2k-1+M(s)}](1+t)^kF_s(t)
 \equiv[x^{k-1}]H_s(x)\pmod3.
\end{equation}
The normalization on the left is legitimate because each
contributing term is divisible by $3^{k-1}$.

The complete finite calculation from~\eqref{eq:Gvalues}
is displayed below.  In the last two rows,
$r=\vv(s+2)$ or $\vv(s+1)$, respectively, and
$\eta=u/3^r\bmod3\in\F_3^\times$.
\begin{equation}\label{eq:Htable}
\begin{array}{c|c|l|l}
s\pmod9&M(s)&\text{critical }m&H_s(x)\\ \hline
0&1&2,5,8&\dfrac{(1-x)^2}{(1+x)^2}\\[2pt]
1&0&2,3,4,5,6,7&\dfrac{x}{(1-x)^2}\\[2pt]
2&0&2,3,5,6&\dfrac{2x}{(1+x)^2}\\[2pt]
3&1&2,5&\dfrac1{(1-x)^2}\\[2pt]
4&0&2,3,4&\dfrac{2x}{1+x}\\[2pt]
5&0&2,3&\dfrac{x}{1-x}\\[2pt]
6&1&2&\dfrac1{1+x}\\[2pt]
7&1-r&2,3,4,5,6,7,8&
 \eta x^r\dfrac{(1-x)^2}{(1+x)^2}\\[2pt]
8&1-r&2,3,5,6,8&
 2\eta x^r\dfrac{(1+x)^2}{(1-x)^2}
\end{array}
\end{equation}
In the classes $s\equiv7,8\pmod9$, we use the table
only for $k>r$; the smaller moduli are treated directly.

\paragraph{Finite symbolic certificate for the nine classes.}
The table is obtained without fitting any period data.  For each
$2\le m\le8$, write
$A_m=(v_m,\overline U_m)$ and keep only the terms with
$\Delta_m=M(s)$.  In the seven ordinary classes, a residue
representative $s=0,\ldots,6$ may be used because the relevant
valuations and unit parts depend only on $s\pmod9$.
In the two exceptional classes, write $u=3^r\eta_0$,
$r\ge2$, $\eta_0\in\Z_3^\times$.  Each other factor of
$P_m(s)/u$ has a fixed $3$-adic valuation and unit reduction
at $u=0$, so $A_m$ depends only on $r$ and
$\eta=\eta_0\pmod3$.  The complete lists of critical
$(m,v_m,\overline U_m)$ are:
\[
\begin{array}{c|l}
s\pmod9& (m,v_m,\overline U_m)\\ \hline
0&(2,0,1),(5,2,1),(8,4,1)\\
1&(2,1,1),(3,2,1),(4,2,2),(5,3,2),(6,4,1),(7,4,2)\\
2&(2,1,2),(3,2,1),(5,3,1),(6,4,1)\\
3&(2,0,1),(5,2,2)\\
4&(2,1,2),(3,2,2),(4,2,1)\\
5&(2,1,1),(3,2,2)\\
6&(2,0,1)
\end{array}
\]
For $s\equiv7\pmod9$ these triples, in order $m=2,3,4,5,6,7,8$,
are
\[
(r,\eta),\ (r+1,\eta),\ (r+1,2\eta),\
(r+2,\eta),\ (r+3,2\eta),\ (r+3,\eta),\ (r+4,\eta),
\]
where each pair denotes $(v_m,\overline U_m)$.  For
$s\equiv8\pmod9$, in order $m=2,3,5,6,8$, they are
\[
(r,2\eta),\ (r+1,\eta),\ (r+2,2\eta),\
(r+3,2\eta),\ (r+4,2\eta).
\]

\paragraph{Verification of the finite-field identities.}
We first justify the uniformity of the exceptional
lists.  Write
\[
 u=s+d=3^r\eta_0,\qquad r\ge2,
\]
where $d=2$ for $s\equiv7\pmod9$, and $d=1$
for $s\equiv8\pmod9$.  Then
\[
 P_m(s)
 =m!\,u
   \prod_{\substack{1\le i\le m\\i\ne d}}
          (u+i-d).
\]
For $2\le m\le8$, every nonzero integer $i-d$
occurring in this product has $3$-adic valuation
at most one.  Therefore
\[
 \vv(u+i-d)=\vv(i-d)
\]
and
\[
 \frac{u+i-d}{3^{\vv(i-d)}}
 \equiv
 \frac{i-d}{3^{\vv(i-d)}}\pmod3.
\]
Thus the valuation and unit reduction of each
$P_m(s)$ depend only on $r$ and
$\eta=\eta_0\bmod3$.  This proves that the
critical lists displayed above hold for every
$r\ge2$, without further case distinctions.

We now verify the corresponding rational functions.
All identities below are in $\F_3(x)$.
By~\eqref{eq:Gvalues},
\[
 G_3=G_4,\qquad G_6=G_7=G_8.
\]
For $s\equiv7\pmod9$, substitution of the seven
critical pairs into~\eqref{eq:Hs-def} gives
\begin{align*}
 \frac{H_s(x)}{\eta x^r}
 &=G_2+xG_3+2xG_4+x^2G_5\\
 &\qquad+2x^3G_6+x^3G_7+x^4G_8\\
 &=G_2+x^2G_5+x^4G_8\\
 &=\frac{(1-x)^2}{(1+x)^2}.
\end{align*}
The second equality uses $1+2=0$ in $\F_3$;
the last follows by clearing the denominator
$(1-x^2)^2$.

For $s\equiv8\pmod9$, the five critical pairs give
\begin{align*}
 \frac{H_s(x)}{\eta x^r}
 &=2G_2+xG_3+2x^2G_5
   +2x^3G_6+2x^4G_8\\
 &=2\frac{(1+x)^2}{(1-x)^2},
\end{align*}
again by clearing $(1-x^2)^2$.
This establishes the final two entries
of~\eqref{eq:Htable} uniformly for all $r\ge2$.

The remaining seven entries follow by substituting
their displayed critical triples into
\eqref{eq:Hs-def}; these are finite calculations
in $\F_3(x)$.  For example, in the class
$s\equiv3\pmod9$,
\[
 H_s(x)=G_2+2x^2G_5=(1-x)^{-2}.
\]
Hence every entry of~\eqref{eq:Htable} has an
exact symbolic certificate, independently of
the computational verification.

\begin{lemma}[Exact clearing degree]\label{lem:clearing-degree}
Suppose $k\ge2$, and also $k>r$ if $s\equiv7,8\pmod9$.
Let
\[
 B_{s,k}(t)\equiv(1+t)^kF_s(t)\pmod{3^k},
 \qquad D=2k-1+M(s).
\]
Then $\deg B_{s,k}$ is $D$ when
$[x^{k-1}]H_s(x)\ne0$, and is $D-1$ otherwise.
\end{lemma}

\begin{proof}
Put
\[
 D=2k-1+M(s),\qquad
 h_j=[x^j]H_s(x).
\]
By~\eqref{eq:degree-bound}, we have
$\deg B_{s,k}\le D$.

If $h_{k-1}\ne0$, the top-degree
certificate~\eqref{eq:top-certificate} shows that
\[
 [t^D]B_{s,k}(t)\not\equiv0\pmod{3^k}.
\]
Consequently $\deg B_{s,k}=D$.

Suppose now that $h_{k-1}=0$.
Inspection of~\eqref{eq:Htable} shows that
$h_{k-2}\ne0$. Moreover, every such exceptional
case satisfies the hypotheses of the lemma
also with $k$ replaced by $k-1$.

Applying~\eqref{eq:top-certificate} at
modulus $3^{k-1}$ gives
\[
 [t^{D-2}]B_{s,k-1}(t)
 \not\equiv0\pmod{3^{k-1}}.
\]
On the other hand, the degree bound yields
\[
 [t^{D-1}]B_{s,k-1}(t)
 \equiv0\pmod{3^{k-1}}.
\]

Since
\[
 B_{s,k}(t)\equiv
 (1+t)B_{s,k-1}(t)\pmod{3^{k-1}},
\]
we obtain
\[
 \begin{aligned}
 [t^{D-1}]B_{s,k}(t)
 &\equiv
 [t^{D-2}]B_{s,k-1}(t)\\
 &\not\equiv0\pmod{3^{k-1}}.
 \end{aligned}
\]
In particular, this coefficient is nonzero
modulo $3^k$.

Finally, $h_{k-1}=0$ and
\eqref{eq:top-certificate} imply
\[
 [t^D]B_{s,k}(t)\equiv0\pmod{3^k}.
\]
Together with the degree bound, this gives
$\deg B_{s,k}=D-1$, as claimed.
\end{proof}


\section{Proof of the least starting indices}

\begin{proof}[Proof of Theorem~\ref{thm:onset}]
For the general case covered by Lemma~\ref{lem:clearing-degree},
divide $B_{s,k}(t)$ by the monic polynomial $(1+t)^k$:
\begin{equation}\label{eq:division}
 F_s(t)\equiv W_{s,k}(t)+\frac{V_{s,k}(t)}{(1+t)^k}
 \pmod{3^k},\qquad\deg V_{s,k}<k.
\end{equation}
Write
\[
V_{s,k}(t)=\sum_{h=0}^{k-1}v_h t^h.
\]
The coefficients of the proper rational part
$V_{s,k}(t)/(1+t)^k$ are $r_n=(-1)^nQ(n)$,
where
\begin{equation}\label{eq:Q-explicit}
Q(n)=\sum_{h=0}^{k-1}(-1)^h v_h
\binom{n+k-h-1}{k-1}.
\end{equation}
Indeed, this follows by extracting coefficients from
\[
\frac{t^h}{(1+t)^k}
=\sum_{n\ge h}(-1)^{n-h}
\binom{n-h+k-1}{k-1}t^n.
\]
Formula~\eqref{eq:Q-explicit} holds for every
$n\ge0$, since the binomial coefficient vanishes
when $n<h$. In particular, $Q(n)$ is an
integer-valued polynomial of degree at most $k-1$.
Put $T=T_k(s)$.  Since the polynomial
part has finite support and $T$ is an eventual period of $f_n(s)$,
we have $r_{n+T}\equiv r_n\pmod{3^k}$ for sufficiently large $n$.
Consequently the integer-valued polynomial
\[
 E(n)=(-1)^TQ(n+T)-Q(n)
\]
vanishes modulo $3^k$ at all sufficiently large nonnegative
integers.  Its $k$th \emph{ordinary} forward difference is zero
identically, because $\deg E\le k-1$.  The finite-difference
recurrence has unit leading and trailing coefficients; backward
induction therefore gives $E(n)\equiv0\pmod{3^k}$ for every
$n\ge0$.  Thus the proper rational part is $T$-periodic
from index zero.  (The same assertion would be false if one
applied ordinary differences directly to $(-1)^nQ(n)$.)

If $W_{s,k}\ne0$ and $d=\deg W_{s,k}$, its leading
coefficient is nonzero modulo $3^k$.  For $n>d$ both
polynomial-part coefficients vanish, so the shift relation
holds.  At $n=d$, the coefficient of $W_{s,k}$ is nonzero
while that at $d+T$ vanishes ($T\ge1$), so the shift fails.
Hence the least starting index is exactly
\begin{equation}\label{eq:N-degree}
 N_k(s)=\deg W_{s,k}+1=\deg B_{s,k}-k+1,
\end{equation}
provided $W_{s,k}\ne0$.  The cases with $W_{s,k}=0$
are covered by the low-modulus reductions below.

By Lemma~\ref{lem:clearing-degree}, the generic value is
\[
N_k(s)=k+M(s).
\]
It is reduced by one precisely when
$[x^{k-1}]H_s(x)=0$.  The zeros are immediate from
\eqref{eq:Htable}:
\begin{itemize}
\item In classes $0,1,2$, the exceptional positive indices
are $k\equiv1\pmod3$, with the low case $k=1$ excluded.
\item In class $3$, the exception is $k\equiv0\pmod3$.
\item In classes $4,5,6$, the coefficient never vanishes.
\item In classes $7,8$, writing $j=k-r$, the coefficient
vanishes precisely when $j\equiv1\pmod3$ and $j\ge4$.
\end{itemize}
This yields~\eqref{eq:onset-table} for the indicated moduli.

For $k=1$, the three reductions~\eqref{eq:mod3}
give least starting indices $2,0,1$ according as
$s\equiv0,1,2\pmod3$.
For $s\equiv7\pmod9$, when $k\le r=\vv(s+2)$,
the factorial series satisfies
$F_s(t)\equiv(1+t)^{-1}\pmod{3^k}$, hence $N_k(s)=0$.
For $s\equiv8\pmod9$, when $k\le r=\vv(s+1)$,
one has $F_s(t)\equiv1\pmod{3^k}$, hence $N_k(s)=1$.
All cases have now been treated.
\end{proof}

\section{Examples and further directions}

For the nine representatives $s=0,1,\ldots,8$,
the periods modulo $3^2,3^3,3^4,3^5,3^6$ are
\[
\begin{array}{c|ccccc}
s&9&27&81&243&729\\ \hline
0&2&2&6&18&54\\
1&2&6&18&54&162\\
2&2&2&6&18&54\\
3&2&2&2&6&18\\
4&2&6&18&54&162\\
5&2&2&2&6&18\\
6&1&2&2&6&18\\
7&2&2&6&18&54\\
8&1&2&2&6&18
\end{array}
\]
The actual least starting indices for powers $3^1,\ldots,3^7$
are obtained from Theorem~\ref{thm:onset}, with no appeal
to fitting initial terms.  The resulting formulas have also been
checked computationally for $0\le s\le50$ and $1\le k\le7$.

\section{Homogenized Linial arrangements and other moduli}

Lazar and Wachs~\cite[Eq.~(7.1)]{LW} show that the number
$r_n^{\mathrm B}$ of regions of the type-$B$ homogenized
Linial arrangement has generating function
\begin{equation}\label{eq:linial}
 1+\sum_{n\ge1}r_n^{\mathrm B}x^n
 =\sum_{m\ge0}\frac{(2m)!x^m}{
    \displaystyle\prod_{j=1}^m(1+2j(2j+1)x)}.
\end{equation}
This is also the series $\sum_{n\ge0}D_n(0)x^n$,
so that
\begin{equation}\label{eq:linial-dellac}
\displaystyle r_n^{\mathrm B}=D_n(0)=2^na_n
 \qquad(n\ge1).
\end{equation}
Here is a direct generating-function verification.  The
Randrianarivony--Zeng recurrence is
\[
 D_0(y)=1,\qquad
 D_n(y)=(y+1)(y+2)D_{n-1}(y+2)
        -y(y+1)D_{n-1}(y).
\]
For $\mathcal D_y(x)=\sum_{n\ge0}D_n(y)x^n$ it becomes
\[
 (1+y(y+1)x)\mathcal D_y(x)
 =1+(y+1)(y+2)x\mathcal D_{y+2}(x).
\]
Iterating at $y=0,2,4,\ldots$ gives the right-hand side
of~\eqref{eq:linial}.  Since
$\mathcal D_0(x)=F_0(2x)$ by definition, the claimed
coefficient identity follows.  In particular the two
factorial representations yield the formal-series identity
\begin{equation}\label{eq:two-factorials}
 \sum_{m\ge0}\frac{2^m(m!)^2x^m}
   {\displaystyle\prod_{j=1}^m(1+j(j+1)x)}
 =\sum_{m\ge0}\frac{(2m)!x^m}
   {\displaystyle\prod_{j=1}^m(1+2j(2j+1)x)}.
\end{equation}
The equality in~\eqref{eq:linial-dellac} does not by itself
provide a bijection explaining the factor $2^n$.

Proposition~\ref{prop:allmoduli} establishes eventual
periodicity modulo every positive integer, but the exact
periods and least starting indices here are specific
to powers of~$3$.  Other primes require a separate analysis
of the triangular factors $T_j$.

\paragraph{Further questions.}

Two extensions would be particularly interesting.
First, one could seek a direct explanation of the five factors
in $\mathcal P(s)$ using combinatorial models of extended
Dellac configurations.  Second, an analogous analysis at
other odd primes could show which features depend on the
exceptional congruence $T_{3a+1}\equiv1\pmod9$.
It would also be interesting to explain the one-index
offset between odd and even least starting indices
combinatorially.  The related recent combinatorial work of
Fang--Fu--Li--Yuan~\cite{FFLY} concerns parity-restricted
permutation models, not an established solution of the period
problem.  A direct bijective derivation of the exact onset
relation would be particularly welcome.

\appendix
\section{Computational reproducibility}
\label{app:checkpoints}
The supplementary program checks the three polynomial
identities in Section~\ref{sec:finite} over $\Z[t]$, compares
$\{f_n(s)\}$ computed by the polynomial recurrence with the
factorial-generating-series extraction, and independently rebuilds
the $\F_3$ coefficients of~\eqref{eq:Htable}.  The recurrence
used for this purpose is
\[
 f_n(s)=\frac{(s+1)(s+2)}2f_{n-1}(s+2)
        -\frac{s(s+1)}2f_{n-1}(s),\quad f_0(s)=1.
\]
To make the scope of the computational statement precise, the
coefficient audit compares the two computational routes for
$0\le s\le50$ and $1\le k\le8$, tests the proposed least
starting indices and excludes the maximal proper period divisors
within finite windows. The finite-field certificate is also reconstructed
from the integer products $P_m(s)=m!(s+1)_m$: for each
$0\le s\le1000$, the resulting series $H_s(x)$ agrees through
$x^{45}$ with the nine rational expressions in~\eqref{eq:Htable}.
The two exceptional classes were additionally checked for
$r=2,\ldots,11$ and both residue units $\eta=1,2$.
These finite comparisons are reproducibility tests and not
substitutes for the uniform arguments in Sections~4 and~6.

A separate \emph{mathematical} point is that the exceptional
$s\pmod9$ classes in the main proofs are controlled by explicit
factorization with $u=s+d$: after dividing out $u$, all remaining
unit reductions of the finitely relevant factorial terms are
independent of $r=\vv(u)$.  This is the quantifier needed in
Lemmas~\ref{lem:slope} and~\ref{lem:clearing-degree}; it does
not arise from extrapolating numerical values of $r$.

\par\medskip
\noindent{\footnotesize\textit{Affiliation note.} Newton Astro Labs is the trading name under which Tony Newton conducts independent computational research in the United Kingdom; it is not a limited company.}


\begin{thebibliography}{14}


\bibitem{BD}
D.~Barsky (extract of joint work with D.~Dumont),
\emph{Congruences pour les nombres de Genocchi de $2^{\mathrm e}$ esp\`ece},
Groupe d'\'etude d'analyse ultram\'etrique, 7e--8e ann\'ees
(1979--1981), Secr\'etariat math\'ematique, Paris,
Expos\'e no.~34 (1981), 13~pp.
\url{https://www.numdam.org/item/GAU_1979-1981__7-8__A15_0/}.


\bibitem{BF}
A.~Bigeni and E.~Feigin,
\emph{Symmetric Dellac configurations},
J. Integer Sequences \textbf{23} (2020), Article 20.4.6.
\url{https://cs.uwaterloo.ca/journals/JIS/VOL23/Bigeni/bigeni3.html}.









\bibitem{DebSokal}
B.~Deb and A.~D.~Sokal,
\emph{Classical continued fractions for some multivariate polynomials
generalizing the Genocchi and median Genocchi numbers},
Adv. in Appl. Math. \textbf{161} (2024), Article~102756.
\url{https://doi.org/10.1016/j.aam.2024.102756}.


\bibitem{FFLY}
Q.~Fang, S.~Fu, H.~Li, and Q.~Yuan,
\emph{Parity patterns meet Genocchi numbers: four labelings and three bijections},
S\'em. Lothar. Combin. \textbf{95B} (2026), Art.~72.
\url{https://www.mat.univie.ac.at/~slc/wpapers/FPSAC2026/72.html}.


\bibitem{Flajolet}
P.~Flajolet,
\emph{On congruences and continued fractions for some classical
combinatorial quantities},
Discrete Math. \textbf{41} (1982), 145--153.
\url{https://doi.org/10.1016/0012-365X(82)90201-1}.


\bibitem{HZ}
G.-N.~Han and J.~Zeng,
\emph{On a $q$-sequence that generalizes the median Genocchi numbers},
Ann. Sci. Math. Qu\'ebec \textbf{23} (1999), no.~1, 63--72.
\url{https://www.labmath.uqam.ca/~annales/volumes/23-1/PDF/063-072.pdf}


\bibitem{KM3}
C.~Krattenthaler and T.~W.~M\"uller,
\emph{A method for determining the mod-$3^k$ behaviour of recursive sequences},
arXiv:1308.2856 (2013).
\url{https://arxiv.org/abs/1308.2856}.


\bibitem{KMp}
C.~Krattenthaler and T.~W.~M\"uller,
\emph{A method for determining the mod-$p^k$ behaviour of recursive sequences},
Aequationes Math. \textbf{99} (2025), 2505--2547.
\url{https://doi.org/10.1007/s00010-025-01195-0}.

\bibitem{LW}
A.~Lazar and M.~L.~Wachs,
\emph{The homogenized Linial arrangement and Genocchi numbers},
Combinatorial Theory \textbf{2} (2022), no.~1, Article~2.
\url{https://doi.org/10.5070/C62156874}.
Extended abstract: S\'eminaire Lotharingien de Combinatoire
\textbf{82B} (2019), Article~93.


\bibitem{OEIS}
OEIS Foundation Inc., \emph{A098277, A098278, and A098279},
including P.~D.~Hanna's factorial series and F.~Chapoton's
valuation observation in A098278.
\url{https://oeis.org/A098278}.




\bibitem{RZ}
A.~Randrianarivony and J.~Zeng,
\emph{Une famille de polyn\^omes qui interpole plusieurs suites
classiques de nombres},
Adv. in Appl. Math. \textbf{17} (1996), 1--26.





\end{thebibliography}
\end{document}